\documentclass[11pt]{article}
\pdfoutput=1 
\usepackage{amssymb}
\usepackage[margin=1in]{geometry}

\usepackage{lineno}
\usepackage{graphicx}
\usepackage{moreverb,url}
\usepackage[ruled,vlined]{algorithm2e}
\usepackage[colorlinks,bookmarksopen,bookmarksnumbered,citecolor=red,urlcolor=red]{hyperref}

\usepackage{xcolor,tabularx,nccmath}
\usepackage{amsmath,amsthm,enumitem}
\usepackage{booktabs,caption,subcaption,mathtools,multirow}

\DeclareMathAlphabet{\mathdutchcal}{U}{dutchcal}{m}{n}
\SetMathAlphabet{\mathdutchcal}{bold}{U}{dutchcal}{b}{n}
\DeclareMathAlphabet{\mathdutchbcal}{U}{dutchcal}{b}{n}

\newcommand{\sB}{\mathdutchcal{B}}

\newcommand{\sK}{\mathdutchcal{K}}

\newcommand{\sT}{\mathdutchcal{T}}

\newcommand{\sR}{\mathdutchcal{R}}
\newcommand{\sC}{\mathdutchcal{C}}
\newcommand{\sS}{\mathdutchcal{S}}

\newcommand{\sZ}{\mathdutchcal{Z}}
\newcommand{\sO}{\mathdutchcal{O}}

\usepackage[nameinlink,capitalise]{cleveref}
\hypersetup{colorlinks={true},allcolors={blue}}

\crefformat{section}{\S#2#1#3}
\crefformat{subsection}{\S#2#1#3}
\crefformat{equation}{#2Equation~\textup{(#1)}#3}
\crefformat{cons}{#2Constraint~\textup{(#1)}#3}
\labelcrefformat{cons}{\textup{#2(#1)#3}}
\crefformat{consset}{#2Constraints~\textup{(#1)}#3} 

\labelcrefformat{consset}{\textup{#2(#1)#3}}
\crefformat{func}{#2Function~\textup{(#1)}#3}
\labelcrefformat{func}{#2Function~\textup{(#1)}#3}
\crefformat{algo}{#2Algorithm~\textup{#1}#3}
\crefformat{objfunc}{#2Objective function~\textup{(#1)}#3}
\labelcrefformat{objfunc}{\textup{#2(#1)#3}}
\crefformat{ch}{#2Chapter~\textup{#1}#3}
\crefformat{sec}{#2Section~\textup{#1}#3}
\crefformat{tab}{#2Table~\textup{#1}#3}
\crefformat{fig}{#2Figure~\textup{#1}#3}
\crefdefaultlabelformat{#2#1#3} 
\Crefname{figure}{Figure}{Figures}
\crefformat{eq}{#2equation~\textup{(#1)}#3}
\labelcrefformat{eq}{\textup{(#2#1#3)}}
\crefformat{theo}{#2Theorem~\textup{#1}#3}
\crefformat{lem}{#2Lemma~\textup{#1}#3}
\crefformat{def}{#2Definition~\textup{#1}#3}
\crefformat{app}{#2\textup{\!#1}#3}
\crefformat{scenario}{#2Scenario~\textup{#1}#3}
 \AtBeginDocument{%
    \crefname{figure}{Figure}{figures}%
}
\let\origref\cref
\def\cref#1{\origref{#1}}

\expandafter\def\expandafter\UrlBreaks\expandafter{\UrlBreaks
  \do\a\do\b\do\c\do\d\do\e\do\f\do\g\do\h\do\i\do\j%
  \do\k\do\l\do\m\do\n\do\o\do\p\do\q\do\r\do\s\do\t%
  \do\u\do\v\do\w\do\x\do\y\do\z\do\A\do\B\do\C\do\D%
  \do\E\do\F\do\G\do\H\do\I\do\J\do\K\do\L\do\M\do\N%
  \do\O\do\P\do\Q\do\R\do\S\do\T\do\U\do\V\do\W\do\X%
  \do\Y\do\Z}

\usepackage{stackengine}

\usepackage{algpseudocode}

\usepackage[normalem]{ulem} 

\newcommand{\replace}[2]{{\color{red}\color{black}{\color{black}#2\color{black}}}} 

\usepackage{rotating}

\crefname{defi}{definition}{definitions}
\Crefname{defi}{Definition}{Definitions}
\crefname{lemma}{lemma}{lemmas}
\Crefname{lemma}{Lemma}{Lemmas}
\crefname{assumption}{assumption}{assumptions}
\Crefname{assumption}{Assumption}{Assumptions}

\usepackage{multicol}

\usepackage{authblk}
\providecommand{\keywords}[1]{\noindent\textbf{Keywords:} #1}

\usepackage[numbers,sort&compress]{natbib}

\begin{document}

\title{Modeling and Calibration of Supplier Selection Problem in Freight Agent-Based Simulations}


\author[1]{Abdelrahman Ismael}
\author[1,2]{Taner Cokyasar}

\affil[1]{Argonne National Laboratory, 9700 S. Cass Avenue, Lemont, IL, 60439 USA}
\affil[2]{Texas A\&M University, College Station, TX 77843, USA}

\maketitle


\begin{abstract}
Freight transportation modeling often struggles with data limitations, especially in accurately representing complex supplier selection processes and their impact on network flows. This research addresses this critical gap by developing a large-scale, calibrated agent-based model for supplier selection, complemented by a probabilistic heuristic for international shipments. Our approach integrates trade relationships between industry sectors, transportation costs, and supplier rating model adapted from existing literature. The model's core objective is to minimize the discrepancy between modeled and observed commodity flows while ensuring a close match to regional shipping distance distributions. Implemented and tested across four major U.S. metropolitan areas, Atlanta, Chicago, Dallas-Fort Worth, and Los Angeles, the model demonstrates high fidelity in replicating observed freight patterns. Key findings reveal consistent alignment with national shipping distance trends and highlight significant spatial variations in commodity trade assignments and demand across the study regions. This behaviorally informed and transport-sensitive framework is designed to approximate real-world decision-making, providing a robust tool for policymakers and planners to evaluate targeted interventions, assess infrastructure investments, and enhance supply chain resilience in the face of disruptions.
\end{abstract}


\section{Introduction}

Freight transportation is a critical enabler of global commerce, facilitating the large-scale movement of raw materials, intermediate goods, and finished products across complex supply chains. These supply chains link producers, distributors, and consumers through multimodal networks that span local, national, and international scale. Efficient navigation of such complexities in global and domestic supply chains is vital to ensure a smooth flow of goods without disrupting the economy. In the United States, for instance, an estimated 20.2 billion tons of freight, valued at over \$18 trillion, moved across the transportation network in 2023—equivalent to roughly 55.5 million tons per day \cite{bts2024}. However, these movements also generate significant negative externalities, including traffic congestion and infrastructure degradation \cite{HOLGUIN2018, Winebrake01082008}. The transportation sector is a critical backbone of the U.S. economy, accounting for a significant percentage of the Gross Domestic Product (GDP) and total logistics costs. Heavy-duty trucks, in particular, serve as the primary mode of freight movement, transporting the vast majority of the nation's shipment value and tonnage \cite{bts_freight_facts}.

In recent years, the vulnerability of global supply chains, particularly to disruptions in freight transportation systems, has become increasingly apparent. High-profile incidents have highlighted how bottlenecks can severely affect the movement of goods and economic stability. The COVID-19 pandemic, for instance, exposed structural weaknesses in logistics systems, as lockdown and labor shortages led to port congestion, vessel delays, and container imbalances, ultimately resulting in widespread supply shortages. Similarly, the 2021 blockage of the Suez Canal—a vital artery for global maritime trade—halted the passage of hundreds of ships and was estimated to delay approximately \$400 million worth of cargo for each hour of the obstruction \cite{larocco_2021, tran_2024}. These events underscore the strategic significance of freight transportation and the necessity for advanced modeling tools that can simulate supply chain dynamics and anticipate the ripple effects of potential disruptions.

Freight transportation is inherently a derived demand, driven by the need to move goods from production to consumption locations. This characteristic highlights the importance of understanding the structure of supply chains and the flow of goods from their origins to destinations, including the logistics and routing choices in between. Changes in production sources can significantly impact transportation networks by altering logistics flows, modal shares, and traffic composition. Conversely, the condition and capacity of transportation infrastructure can influence sourcing decisions, for instance, port congestion or deteriorating network performance may discourage sourcing from specific regions. Additionally, the push for supply chain resilience and nearshoring strategies has begun to reshape global sourcing behavior. For instance, evolving trade dynamics have recently altered U.S. import patterns. In 2023, Mexico overtook China as the United States’ top trading partner \cite{swanson_romero_2024}, signaling a structural shift in freight origins These shifts carry significant transportation implications. Nearly 40\% of Chinese imports to the U.S. enter via the Ports of Los Angeles and Long Beach \cite{tonlexing_2025}, whereas most imports from Mexico enter through Texas by truck and rail \cite{bts_2017}. As sourcing patterns evolve, they necessitate reevaluation of infrastructure investments, highlighting the need for freight modeling tools to accurately assess network-level impacts.

Effectively addressing the network-level impacts of freight movement requires modeling tools capable of capturing the decision-making behavior of individual freight agents, particularly with regard to sourcing and logistics choices. Recently, agent-based modeling (ABM) has been utilized in the freight domain as a mean of simulating such complex interactions. In contrast to traditional aggregate or deterministic models, freight ABMs represent suppliers, shippers, carriers, receivers, and end consumers as autonomous agents with unique objectives and behavioral rules. These agents interact with each other and with their physical and policy environment in ways that give rise to emergent system-level outcomes. This approach enables the modeling of heterogeneous preferences, adaptive behaviors, and decentralized decision-making. As a result, freight ABMs are especially well-suited to exploring how logistical decisions—such as supplier selection, carrier choice, shipment size, transport mode, and routing—are shaped by evolving conditions including infrastructure constraints, regulatory policies, technological innovation, and interactions with passenger traffic on shared networks.

While freight ABMs encompass a wide range of logistics behaviors, this study focuses specifically on the sourcing choices of freight agents, namely, supplier selection and commodity assignment. These sourcing choices represent a foundational decision in supply chain operations, shaping not only procurement outcomes but also influencing freight demand and the spatial distribution of goods movement. From a transportation modeling perspective, these choices determine the origins of commodity flows, which in turn affect network usage, traffic patterns, and infrastructure conditions. Opting for suppliers located near demand centers can minimize vehicle miles traveled (VMT), fuel consumption, and associated economic burden. However, in many cases, commodities are available only from geographically distant or cost-advantaged regions, necessitating long-haul shipments and greater reliance on national freight corridors. These trade-offs highlight the importance of explicitly incorporating supplier selection into freight ABMs to better simulate real-world logistics dynamics and evaluate the system-level consequences of sourcing behavior.

This research presents a large-scale calibrated ABM of supplier selection and commodity assignment (within POLARIS simulation tool) that integrates trade relationships between industry sectors and transportation costs. The model aims to emulate a more realistic representation of how receiver businesses choose their suppliers and how commodities flow between each pair of suppliers and receivers, \replace{}{compared to previous models (including earlier versions of POLARIS-Freight)}. This is being achieved by accounting for receiver's own benefit of maximizing their perception of the selected suppliers' rating using a supplier rating model from the literature \cite{pourabdollahi_2017}. Additionally, the model seeks to account for unobserved factors through ensuring that the inter-zonal flows are matched and the shipping distance distribution gap is reduced. By linking micro-level attributes of shipping cost and supplier rating with macro-level freight patterns, the model offers a good representation of freight flows between agents and allow modeling and analyzing impacts of freight transportation on infrastructure planning and supply chain resilience. In this paper, the model is implemented in four metropolitan areas: Atlanta, Chicago, Dallas-Fort Worth (DFW), and Los Angeles (LA) along with a heuristic to select importer and exporter establishments in these four areas.

The rest of this paper is organized as follows. \nameref{literature} section provides details on traditional supplier selection models, supplier selection in freight ABMs, and research gap and contributions. \nameref{framework} section lays out the supplier and commodity selection module within POLARIS Freight and points to relevant data sources that are utilized in such studies. \nameref{methodology} section provides model notations and algorithmic details for both domestic and international trade models. \nameref{results_discuss} section presents findings for four metro areas. Finally, findings, limitations, and policy implications are explained in \nameref{conclusion} section.

\section{Literature Review}\label{literature}
This section describes a brief summary of research work done in freight ABMs and supplier selection modeling in the literature.

\subsection{Traditional Supplier Selection Models}

Supplier selection is a well-described problem in supply chain management and operations research, with a rich body of literature addressing how receivers choose among potential suppliers. Originally supplier selection models used minimum cost as the most important criteria, however, later research started including other attributes such as geographic location \cite{sarkar_mohapatra_2006, thanaraksakul}, quality \cite{chan_kumar_2007}, delivery \cite{sarkar_mohapatra_2006}, performance history \cite{watt_kayis_willey_2010}, reputation \cite{sarkar_mohapatra_2006, thanaraksakul, watt_kayis_willey_2010}, risk \cite{chan_kumar_2007}, service levels \cite{chan_kumar_2007}, production capacity \cite{sarkar_mohapatra_2006}, technology \cite{chan_kumar_2007, thanaraksakul}, etc. These models generally fall into multi-criteria decision-making (MCDM) and mathematical optimization.
MCDM techniques such as the Analytic Hierarchy Process (AHP) \cite{ishizaka_2012, chen_chao_2012}, Technique for Order of Preference by Similarity to Ideal Solution (TOPSIS) \cite{crispim_2008, önüt_2009, azadeh_2010}, Analytic Network Process (ANP) \cite{demirtas_2009, önüt_2009}, Data Envelopment Analysis (DEA) \cite{saen_2010, azadeh_2010}, among others, are widely used to capture subjective preferences and trade-offs.
On the other hand, optimization-based approaches—including single- and multi-objective linear programming \cite{razmi_2009, demirtas_2009}, goal programming \cite{demirtas_2009}, and mixed-integer programming—have been applied to more quantitative supplier selection scenarios. A comprehensive review of these supplier selection models is presented in multiple papers including \cite{de_boer_labro_morlacchi_2001, chai_liu_ngai_2013, taherdoost_brard_2019}.

\subsection{Supplier Selection in Freight Agent-Based Models}

While traditional supplier selection models focus on procurement efficiency, freight ABMs aim to simulate the behavior of suppliers and receivers within transportation systems. In recent years, freight ABMs have been gaining traction and more researchers started developing different ABMs, including FAME  \cite{samimi_2014}, MASS-GT \cite{de_bok_2018, de_bok_2025}, CRISTAL \cite{stinson_mohammadian_2022}, SynthFirm \cite{spurlock_2024}, and POLARIS \cite{auld_2016, zuniga_2023, ismael_2025}. 

The Freight Activity Microsimulation Estimator (FAME) \cite{samimi_2014} models supply chains at an aggregated firm level. FAME aggregates firms based on their locations, type, and size to generate firm-type synthetic population. The model uses a fuzzy-rule based model for supplier selection, where the rules depends on categorical variables for size and proximity of suppliers and receivers. MASS-GT \cite{de_bok_2018, de_bok_2025} follows a top-down approach, where they synthesize shipments from aggregated Origin-Destination (OD) flows, and then assign them to a given receiver followed by a given supplier. The receiver probabilistic assignment is a function of the probability of a shipment being used by a specific industry sector and the firm sizes within this sector. While, the supplier assignment is a function of the probability of a shipment being sent by a specific industry sector and the firm sizes within this sector weighted by the transportation cost. SynthFirm \cite{spurlock_2024} synthesizes firms and uses a market-clearing mechanism based on shipping distances and values estimated from the Commodity Flow Survey (CFS) data. 

CRISTAL \cite{stinson_mohammadian_2022} models long-term, medium-term and short-term horizons to capture supply chain behaviors in-depth. Long term decisions include how much goods to produce and how much goods are needed, fleet and warehousing decisions, and trade partnerships including supplier selection. The medium term decisions include setting up order frequencies and tour generation, while the short term decisions include powertrain and routing choices. The supplier selection model in CRISTAL is a Multinomial Logit Model (MNL) where the utility associated with a given supplier out of a candidate set of suppliers, is based on the following supplier attributes: if it is located in the same region (internal), if it is a foreign supplier, employment, great circle distance (GCD) from the receiver location. The parameters were not estimated based on a collected survey data, their values were based on a fuzzy logic model from a report of \citep{cambridge_systematics_2011}.

POLARIS \cite{auld_2016} is an agent-based and activity-based model for both passengers and freight. The passenger modules of POLARIS synthesize the population and their activities, estimate activity destinations and resolve conflicts between activity start times. POLARIS also estimates mode choice and routing decisions, and simulates passenger cars, transit, ride-hailing services, etc. The POLARIS Freight modules \cite{zuniga_2023, ismael_2025} is mostly an enhanced version of CRISTAL. POLARIS-Freight uses the core modules of CRISTAL and improves the structure of models with updated models and collected data. POLARIS-Freight also adds other freight modules such as passing through freight loaded and empty demand, service trip demand, on-demand deliveries (ODD) for meal and groceries, port allocation, among others. POLARIS-Freight offers the additional benefit of interacting with the passenger demand, so e-commerce and  ODD demands are a function of the synthesized population attributes, and directly impact the shopping and eat-out activities of the households. Moreover, it allows the co-simulation of freight and passenger trips which account for their inter-dependencies and interactions on the traffic network and their en-route switching decisions due to congestion. The current paper proposes a large-scale supplier selection model to replace the utility based model of CRISTAL existing in POLARIS-Freight.

Another important ABM study \cite{pourabdollahi_2017} that strongly relates to this research, combined decision-making framework with computational techniques for supplier selection. The authors introduced a hybrid agent-based computational economics and optimization model, demonstrating how behavioral rules and optimization techniques can be integrated to simulate decentralized procurement decisions. The authors used fuzzy logic and genetic algorithms to explore large solution spaces. This paper is particularly important for our current research, since we adapted the supplier rating utility model developed in \cite{pourabdollahi_2017} based on collected survey data to account for factors such as reliability, financial credibility, and capacity. 

\subsection{Research Gap and Contributions}

The integration of supplier selection into freight ABMs remains an open challenge. Despite their sophistication, traditional supplier selection models often abstract away the behavioral aspects of individual agents when selecting suppliers and the flow calibration which is critical from a transportation systems perspective. \replace{Conversely, supplier selection in most freight ABMs is treated exogenously. Decisions are not typically informed by network conditions, and are not constrained by actual zonal commodity flows.}{Conversely, supplier selection in most freight ABMs, including previous iterations of POLARIS-Freight and CRISTAL, relies on exogenous assignments or unconstrained utility maximization. By introducing strict calibration discipline to macroscopic trade targets and an optimization-based decomposition that scales to millions of candidate pairs, this framework directly removes the limitation of uncalibrated trade flow assignments, ensuring decisions are both informed by network conditions and constrained by actual zonal commodity flows.} Moreover, international trade flows are often ignored or aggregated on a zonal level with little regard to which specific businesses act as importers or exporters and which ports serve as the points of entry. Not all models consider detailed supplier attributes, e.g., reliability, credibility, capacity, when modeling procurement decisions. This limits the ability of freight ABMs to represent a more realistic picture of the transportation component of supply chains, in order to effectively evaluate infrastructure investments, policy shifts, or supply chain resilience strategies.

This research aims to bridge these two domains by introducing a calibrated, large-scale supplier selection and commodity assignment model into a freight ABM context. Specifically, the primary objectives of this study are to i) move beyond exogenous supplier assignment by developing an optimization-based framework that accounts for firm-level behaviors, including reliability ratings and capacity constraints,
ii) ensure that the modeled supplier-receiver pairings reproduce observed macroscopic patterns, specifically regional shipping distance distributions and zonal commodity flow volumes, iii) overcome the limitations of zonal aggregation by establishing a heuristic to assign import/export flows to specific domestic establishments based on port throughput and industry sector data. Our contribution is threefold:
\begin{enumerate}
    \item \textit{Optimization-Based Supplier Selection and Commodity Assignment}: We introduce a linear programming formulation that selects supplier–receiver pairs by commodity for domestic trade based on shipping cost, supplier production capacity, receiver consumption needs, commodity compatibility, and supplier rating, adapted from \cite{pourabdollahi_2017}, which model indirectly cost of commodity, reliability, and financial credibility.
    \item \textit{Heuristic International Assignment}: We develop a probabilistic heuristic that selects importers and exporters for international shipments, grounded by individual port flows and North American Industry Classification System (NAICS) industry sector to commodity mappings.
    \item \textit{Model Calibration and Large-scale Implementation}: The model components are tested on four metropolitan areas in the U.S: Atlanta, Chicago, DFW, and LA, and calibrated to match observed trade patterns and commodity flows. The model is integrated within POLARIS, yielding additional benefits from utilizing synthesized freight agent attributes.
\end{enumerate}

\section{Research Framework and Data Sources}\label{framework}

\begin{figure}
    \centering
    \includegraphics[width=\linewidth]{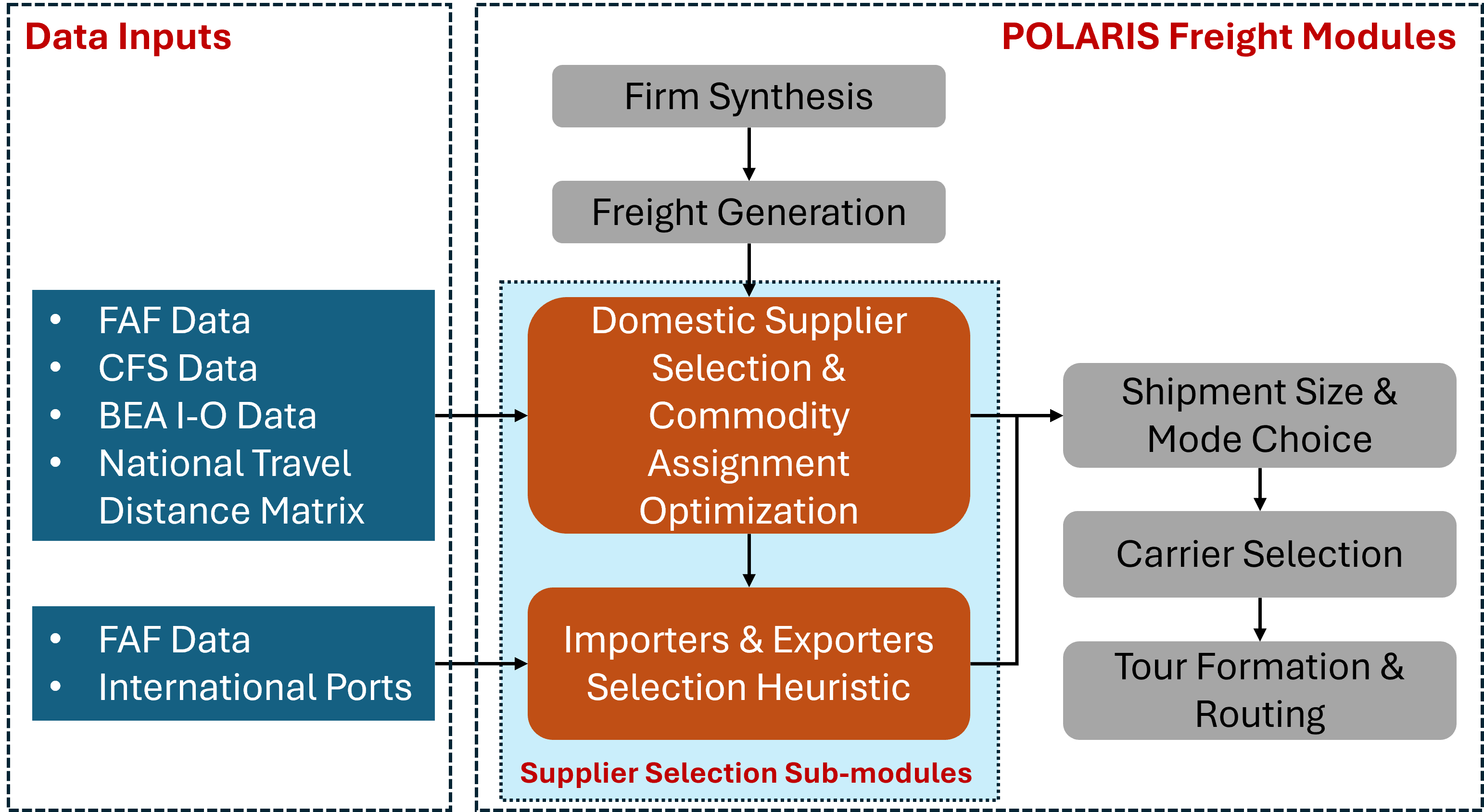}
    \caption{Supplier and commodity selection module within POLARIS Freight.}\label[fig]{sup_sel_framework}
\end{figure}

The data inputs for this model include various public data sources and POLARIS model outputs, as shown in \cref{sup_sel_framework}.

\subsection{Public Data Sources}
    \begin{itemize}
    \item The Bureau of Economic Analysis Input–Output data \cite{bea_2025} contains information on make-use interactions between industry sectors. These interactions capture how sectors produce commodities for others, consume commodities produced by other sectors, and the volume of trade that occurs among them. This data source is used in the model to identify potential sectors that can supply commodities to the industry sector of the receiver. \replace{}{In the absence of observed firm-to-firm logistics relationships, these Make-Use tables serve as deterministic structural constraints governing the feasible pairings between firms belonging to different NAICS codes.}
    
    \item Freight Analysis Framework (FAF) data \cite{bts_2017} maps the commodity flow tonnage among 140 FAF zones, including: 132 domestic FAF zones (representing major metropolitan areas by state and rest of their states) and 8 foreign FAF zones (representing countries like Canada and Mexico, sub-continents or whole continents). This data source is used in the model to provide the zonal commodity flows to be matched by the model. The commodities are classified according to the 2-digit Standard Classification of Transported Goods (SCTG) codes.

    \item Commodity Flow Survey (CFS) public use file (PUF) data \cite{CFS} represent a sample from survey data for shippers exporting or shipping commodities domestically including the shipment weight and distances along with weight factors for the sample. This data source is used in the model to provide a shipping distance distribution for each city to be matched by the model. The commodities are classified according to the 2-digit SCTG codes.
\end{itemize}

Given the mixture of usage of NAICS industry sectors and SCTG commodities in the above data, and that businesses, i.e., suppliers and receivers are typically classified using NAICS codes, it is crucial to use a mapping between the two classifications. It is important to mention, however, that there is no standard one NAICS to many commodities mapping that can be used for all businesses with the same NAICS code, since businesses from the same industry sector can produce different mixtures of commodities. Acknowledging this limitation, we opted for using one of the mappings in the literature developed in \cite{pourabdollahi_2015}. In addition, we aggregated the 42 SCTG codes into 15 groups shown in \cref{commodity_groups}, this aggregation helps in reducing the problem size. However, it comes at the expense of capturing the detailed behavior of the supply chains of individual commodities.

\begin{table}[ht]
\centering
\caption{Commodity Grouping Used}
\label[tab]{commodity_groups}
\begin{tabularx}{\linewidth}{lX}
\hline
\textbf{Label} & 
\textbf{Commodity Group Name} \\ 
\hline
1  &	Food, Agriculture, and Forestry Products \\
2  &	Mining Products \\
3  &	Petroleum Products \\
4  &	Chemical and Pharmaceutical Products \\
5  &	Wood Products \\
6  &	Paper Products \\
7  &	Nonmetallic Mineral Products \\
8  &	Metal and Machinery Products \\
9  &	Electronic, Electrical and Precision Equipments \\
10 &	Motorized and Transportation Vehicles and Equipments \\
11 &	Household and Office Furniture \\
12 &	Plastic, Rubber and Miscellaneous Manufactured Products \\
13 &	Textiles and Leather Products \\
14 &	Waste and Scrap \\
15 &	Mixed and Unknown Freight \\
\hline
\end{tabularx}
\end{table}

\subsection{POLARIS Model Outputs}
\begin{figure}
    \centering
    \includegraphics[width=\linewidth]{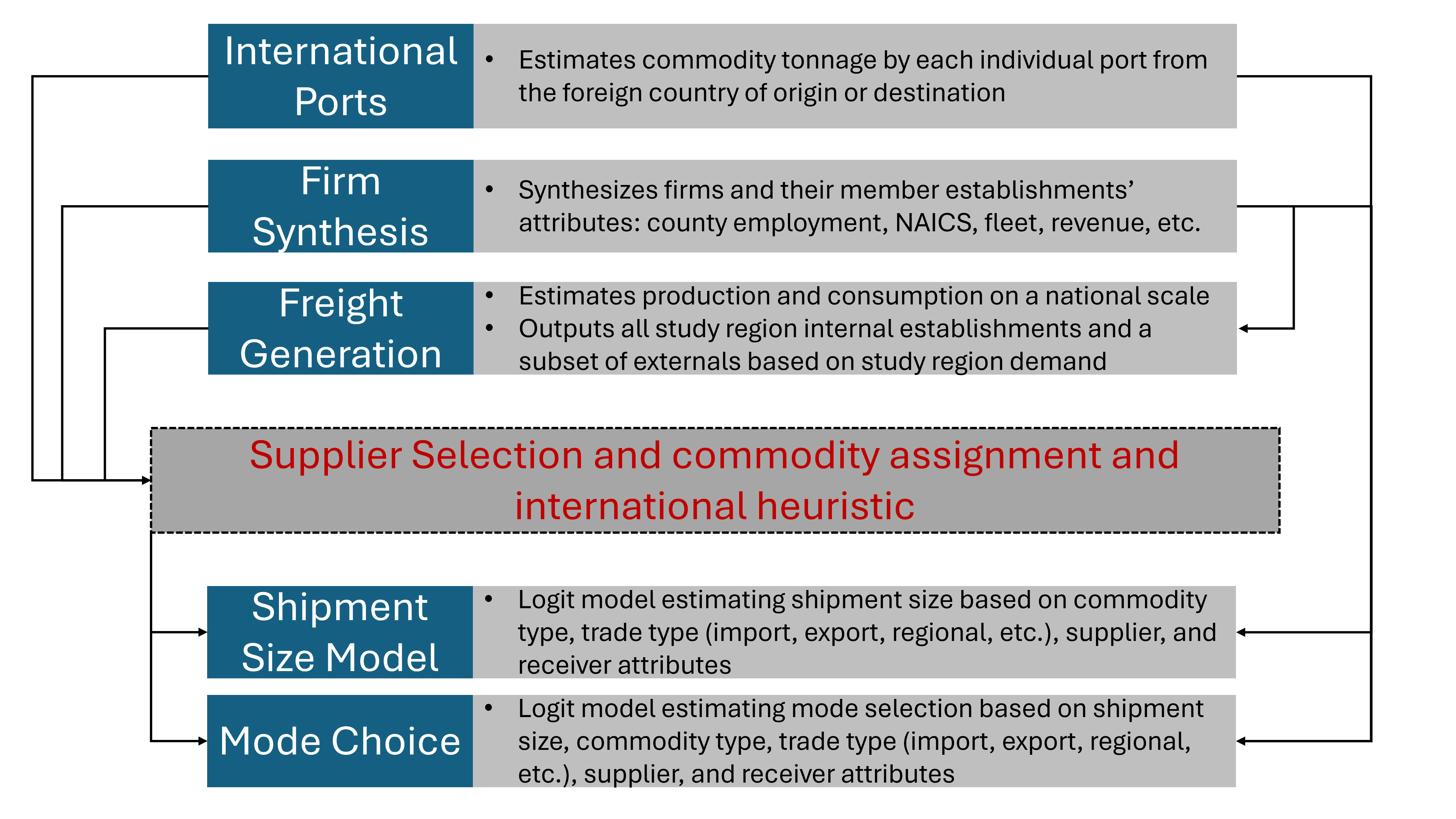}
    \caption{POLARIS Freight modules interactions.}\label[fig]{sup_sel_flowchart}
\end{figure}

While this model is capable of independently addressing the supplier selection problem, implementation within POLARIS can offer significant advantages. As this integration allows access to different attributes for businesses through POLARIS' firm synthesis and freight generation modules, as shown in \cref{sup_sel_flowchart}. This integration can also enable future dynamic feedback between sourcing decisions and network conditions that can be translated into modified shipping costs. 

    \begin{itemize}
    \item \textit{Firm Synthesis}: 
    This module synthesizes parent firms and their member establishments (business locations) along with their characteristics based on proprietary data samples to match control totals from public data sources. The synthesized attributes include 3-digit NAICS industry sectors, U.S. county, employment, fleet, revenue, etc.
    
    \item \textit{Freight Generation (FG)}: 
    This module uses the firm synthesis outputs, FAF data, and NAICS to commodity mapping to estimate the production capabilities and consumption demand of each establishment. The production and consumption of establishments are computed based on tonnage rates per employee, where rates differ based on FAF zones and NAICS.
    
    Given the estimated production and consumption capacities of external establishments, we randomly sample from each external zone enough establishments to cover that zone's total supply and demand to/from the study region. This helps reduce the problem size, where a given receiver will not use all possible national suppliers as potential suppliers. Hence, for a given study region, all internal establishments within the region are considered, however, only a portion of the external domestic establishments are used. This is a trade-off as it reduces the problem size but in the same time limits the selection pool of external suppliers.
    
    \item \textit{International Ports}: 
    For the international shipment heuristic, POLARIS disaggregates import and export FAF flows from the zonal level to individual ports using ports and land borders information and capacities from Bureau of Transportation Statistics \cite{bts_2025}. These ports are used as points of entry and exit for import and export flows for the businesses in the study region.
    \end{itemize}

These POLARIS Freight module interactions are summarized in in \cref{sup_sel_flowchart}, where the POLARIS model outputs feed into the supplier selection and commodity assignment problem, resulting in annual trade flows between suppliers and receivers. These trade flows include information on the quantity of tonnage traded annually, type of commodities traded, and trade type. The trade type depends on the type and location of supplier and receiver: import, export, regional, domestic inbound (external-internal), domestic outbound (internal-external). Shipment size and mode choice logit models depend on such information to estimate shipping chain decisions. 

\section{Methodology}\label{methodology}

In this section, we first begin with formulating a linear programming model to solve the supplier selection and commodity assignment problems jointly. Since commercial solvers struggle with computationally modeling this problem at the target area scales, we decompose the problem into two phases, addressing the supplier selection problem first, and solving the commodity assignment problem next. \cref{notations} provide definitions for sets, parameters, and variables used in this section. To facilitate readability, the model adheres to the following structural logic: Uppercase Roman letters denote parameters (e.g., $C_{sr}$), while Lowercase Roman letters denote variables (e.g., $x_{sro}$), with the exception of the weight parameter $w$. Superscripts are reserved strictly for modifiers defining specific subsets (e.g., $\mathcal{R}^\sigma$). Furthermore, subscripts consistently follow the flow of goods: \textit{Origin $\to$ Destination $\to$ Commodity}. For example, the triplet $(s,r,o)$ always refers to supplier $s$, receiver $r$, and commodity $o$, respectively.

\subsection{Joint Supplier Selection and Commodity Assignment}

Let $\sR$ and $\sS$ denote a set of receivers and a set of suppliers, respectively. Also, let $\sO$ represent a set of commodity types. The goal in the joint supplier selection and commodity assignment problem is to find the percent of commodity type demand $o\in\sO$ of receiver $r\in\sR$ by supplier $s\in\sS$, denoted by $x_{sro}\in [0,1]\land x_{sro}\in\mathbb{R}_{\geq 0}$. While matching as much as receivers and suppliers is a target, we also desire to minimize the percent of unmet demand of receiver $r$, $u_r\in [0,1]\land u_r\in\mathbb{R}_{\geq 0}$. We categorize every potential $(s,r)$ pair into a distance bin $b\in\sB$. Then, we define the percentage of observed tonnage in distance bin $b\in\sB$ from total observed regional tonnage from the CFS data denoted by $Q_b$. The gap in demand tonnage in distance bin $b\in\sB$ is represented by $g_b$, that is the absolute difference of demand of receiver $r\in\sR$ denoted by $D_r$ supplied by all suppliers in distance bin $b\in\sB$ divided by total demand. Another variable we seek the value for is $y_{zz^\prime o}$ defined as the absolute gap in the tonnage flow of commodity $o\in\sO$ from zone $z\in\sZ$ to zone $z^\prime \in\sZ$. The value for this variable depends on the assignments of $x_{sro}$ and $P_{zz^\prime o}$, which is defined as the flow of goods in tonnage of commodity $o\in\sO$ from zone $z\in\sZ$ to zone $z^\prime \in\sZ$. Finally, we account for rating of receiver $r\in\sR$ for supplier $s\in\sS$, that is how likely they would be matched, denoted by $N_{sr}$. We then present the linear program as follows:

\begin{table}[!htb]
  \footnotesize
  \caption{Sets, parameters, and variables used in Supplier Selection and Commodity Assignment models} \label[tab]{notations}
  \begin{tabularx}{\linewidth}{lX}
  \toprule
    \textbf{Set} & \textbf{Definition}\\
    \midrule
    $\sB$ & set of distance bins\\
    $\sK_b$ & subset of supplier $s\in\sS$ and receiver $r\in\sR$ where the distance in between is in the distance bin $b\in\sB$, that is the subset for a given $b\in\sB$ returns $(s,r)$\\
    $\sO$ & set of 2-digits SCTG commodity groupings \\
    $\sO_{sr}$ & set of possible commodities that can be supplied by supplier $s\in\sS$ to receiver $r\in\sR$ \\
    $\sR$ & set of receivers\\
    $\sR_s$ & subset of receivers can be supplied by supplier $s\in\sS$\\
    $\sR^\sigma$ & subset of micro size receivers, $\sR^\sigma \subset \sR$\\
    $\sS$ & set of suppliers\\
    $\sS_r$ & subset of suppliers that can supply receiver $r\in\sR$\\
    $\sZ$ & set of transportation analysis zones (TAZs)\\
    \midrule
    \textbf{Parameter} & \textbf{Definition}\\
    \midrule
    $C_{sr}$ & transportation cost of shipment between supplier $s \in \sS$ and receiver $r \in \sR$ \\
    $D_{r}$ & demand of receiver $r\in\sR$ \\
    $K_{s}$ & supply capacity of supplier $s\in\sS$\\
    $N_{sr}$ & rating of receiver $r\in\sR$ for supplier $s\in\sS$ \\
    {$P_{z z^\prime o}$} & {flow of goods (tonnage) of commodity $o \in \sO$ from zone $z\in\sZ$ to in zone $z^\prime \in \sZ$} \\
    $Q_{b}$  & percentage of observed tonnage in distance bin $b\in\sB$ from total observed regional tonnage from the CFS data \\
    $W_{sr}$ & amount of demand of receiver $r\in\sR$ supplied by supplier $s\in\sS$\\
    $w_1^r$ & objective function weight for unmet demand of receiver $r\in\sR$ \\ 
    $w_2$ & objective function weight for shipping cost between all pairs of supplier $s \in \sS_r$ and receiver $r \in \sR$  \\ 
    $w_3$ & objective function weight for the rating of all receivers $r \in \sR$ for suppliers $s \in \sS_r$\\ 
    $w_4$ &  objective function weight for absolute percentage gaps in all distance bins $b \in \sB$ between modeled and observed shipping distance distribution\\ 
    $w_5$ &  objective function weight for absolute gaps between all modeled and observed commodity flows\\ 
    $Z_i$ & zone in which establishment $i\in \sR \cup \sS$ operates\\
    \midrule
    \textbf{Variable} & \textbf{Definition}\\
    \midrule
    $c_{sro}$ & percent of commodity $o \in \sO$ assigned to the pair of supplier $s\in\sS$ and receiver $r\in\sR$ \\
    $g_{b}$ & percent gap of demand tonnage in a distance bin $b \in\sB$ \\
    $u_{r} $ & percent of unmet demand of receiver $r \in\sR_t$ \\
    $x_{sr}$ & percent of demand of receiver $r\in\sR_t$ met by supplier $s\in\sS$\\
    $x_{sro}$ & percent of commodity type demand $o\in\sO$ of receiver $r\in\sR$ by supplier $s\in\sS$\\
    {$y_{z z^\prime o}$} & {gap in the tonnage flow of commodity $o \in \sO$ from zone $z \in \sZ$ to zone $z^ \prime \in \sZ$} \\
    \bottomrule
  \end{tabularx}
\end{table}

\begin{equation}\label[objfunc]{obj_fun_joint}
\begin{split}
    min_{x,u,g, y} & \quad w_1^r \sum_{r\in\sR, s\in\sS} C_{sr} D_{r}  u_{r} + w_2 \sum_{r \in \sR, s \in \sS} C_{sr} D_{r} \sum_{o \in \sO} x_{sro}  \\ & - w_3 \sum_{r \in \sR, s \in \sS} N_{sr} \sum_{o \in \sO} x_{sro} + w_4  \sum_{b \in \sB} g_{b} \\ & + w_5 \sum_{z \in \sZ, z^ \prime \in \sZ , o \in \sO} y_{z z^\prime o}
\end{split}
\end{equation}

where \[w_1^r = 
\begin{cases}
    10 w_1,      & \text{if } r \in R^\sigma \\
    w_1,    & \text{otherwise}
\end{cases}
\]

\begin{equation}\label[consset]{FM_Supplier}
    \sum_{r \in \sR, o \in \sO} D_{r}x_{sro}  \leq K_{s} \quad \forall s \in \sS 
\end{equation}

\begin{equation}\label[consset]{FM_unmet_demand}
    \sum_{s \in \sS, o \in \sO} x_{sro} + u_{r} = 1  \quad \forall r \in \sR
\end{equation}
\begin{equation}\label[consset]{FM_flow_deviation}
    g_{b} = \left | \frac{\sum_{(s,r)\in\sK_b} D_{r} \sum_{o \in \sO} x_{sro}}{\sum_{r\in\sR} D_r}  - Q_{b} \right |, \quad \forall b \in \sB
\end{equation}

\begin{equation}\label[consset]{FM_FAF_deviation}
    y_{z z^\prime o} = \left | \sum_{\substack{s \in \sS|Z_s=z,\\ r \in \sR|Z_r=z^\prime}} D_r x_{sro} - P_{z z^\prime o}  \right |, \quad \forall z, z^ \prime \in \sZ, o \in \sO
\end{equation}

\begin{equation}\label[consset]{FM_commodity_non_negativity}
    x_{sro}, u_{r}, g_{b}, c_{sro}, y_{z z^\prime o}  \in \mathbb{R}_{\geq 0}, x_{sro}, u_{r}, c_{sro} \in [0,1]
\end{equation}

In objective function \labelcref{obj_fun_joint}, we minimize the weighted cost of i) unmet demand, ii) the transportation cost for met demand, iii) supplier-receiver rating factor, iv) percentage gap in demand tonnage in distance bins, and v) gap in the inter-zonal commodity tonnage flow. It is important to note that the components of objective function possess vastly different orders of magnitude (e.g., total transportation cost versus percentage-based flow gaps). To prevent the larger magnitude terms from mathematically dominating the optimization, the problem is solved using a hierarchical (lexicographic) multi-objective framework. Implemented via the \textit{setObjectiveN} feature in the Gurobi solver \cite{gurobi}, this approach optimizes objectives sequentially based on strict priority levels rather than a simultaneous weighted sum. Priority is assigned in the following order: (1) minimizing unmet demand to ensure system stability, (2) minimizing the deviation from observed shipping distance distributions, and (3) minimizing transportation costs and maximizing supplier ratings. \replace{}{The normative choice to prioritize calibration over cost minimization ensures that the model bounds theoretical micro-economic efficiency within the realities of observed macroeconomic behavior. This sequential structure often creates binding trade-offs where the solver intentionally bypasses a highly-rated, low-cost local supplier if selecting them would cause the aggregate shipping distance distribution to violate the observed regional patterns. Ultimately, this method inherently handles scaling disparities without requiring heuristic normalization.} This sequential solution method inherently handles the scaling disparities without requiring heuristic normalization of the input data. To prioritize micro size receivers denoted by $\sR^\sigma$, we inflate $w_1$ by 10. Note that the number used here could be parameterized, and a different number might be more appropriate for another dataset. Constraints \labelcref{FM_Supplier} ensures that the supply provided by a given supplier $s\in\sS$ does not exceed its supply capacity $K_s$. Constraints \labelcref{FM_unmet_demand} satisfies the integrality of the total demand received by $r\in\sR$. Constraints \labelcref{FM_flow_deviation} define $g_b$ in terms of $x_{sro}$ variables. Constraints \labelcref{FM_FAF_deviation} define $y_{zz^\prime o}$ in terms of $x_{sro}$ variables. Finally, non-negativity constraints \labelcref{FM_commodity_non_negativity} denote variables and their domains.

Collectively, the objective function captures the tension between \textit{theoretical efficiency} and \textit{empirical realism}. The terms related to transportation cost ($C_{sr}$) and supplier rating ($N_{sr}$) represent the micro-economic incentives of individual agents to minimize expenses and maximize service quality. Conversely, the calibration terms ($g_b$) represent the macro-level structural constraints of the freight system. By penalizing deviations from observed distance distributions, the model accounts for unobserved logistical frictions, such as long-term contracts, supply chain inertia, or specialized commodity compatibility, which often prevent agents from selecting the absolute closest supplier. Thus, the optimization seeks the most efficient sourcing configuration that remains consistent with historical trade patterns.

While the model is formulated using percentage flows ($x_{sro} \in [0,1]$) for computational scaling, it is helpful to conceptually relate these to absolute flows for clarity. Let $T_{sro}$ represent the absolute tonnage supplied by supplier $s$ to receiver $r$. The relationship between the two is defined as $T_{sro} = D_r x_{sro}$. Consequently, the standard demand constraint $\sum_{s} T_{sro} = D_r$ is normalized by dividing both sides by the total demand $D_r$, resulting in the utilized constraint $\sum_{s} x_{sro} = 1$. This normalization is critical for the solver's performance, as it bounds variables within a unit scale $[0,1]$, preventing numerical instability caused by the high variance in demand magnitudes ($D_r$) across different receiver types.

\subsubsection{Supplier Selection Model}

The joint model though linear and seemingly simple is not scalable for the problem sizes tackled in this paper due to very large number of potential combinations for $(s,r,o)$ triplet. To this end, we decompose the problem into supplier selection and commodity assignment. The details of the decomposition and more details on how deal with large-scale instances of the problem are presented in \cref{supplier_selection_algo}. \replace{}{It should be noted that the optimization is performed concurrently for establishments across all industry sectors within these spatial subproblems, rather than solved separately by NAICS code. This concurrent approach is essential to capture realistic cross-sector competition for shared supplier capacities.} The portion of the model dealing with the supplier selection is as follows.

\begin{equation}\label[objfunc]{obj_fun_ss}
\begin{split}
    min_{x,u,g} & \quad w_1^r \sum_{r\in\sR, s\in\sS} C_{sr} D_{r}  u_{r} + w_2 \sum_{r \in \sR, s \in \sS} C_{sr} D_{r}  x_{sr} \\ & - w_3 \sum_{r \in \sR, s \in \sS} N_{sr} x_{sr}  + w_4  \sum_{b \in \sB} g_{b}
\end{split}
\end{equation}

where \[w_1^r = 
\begin{cases}
    w_1,      & \text{if } r \in R^\sigma \\
    10w_1,    & \text{otherwise}
\end{cases}
\]

\begin{equation}\label[consset]{Supplier}
    \sum_{r \in \sR} D_{r}x_{sr}  \leq K_{s} \quad \forall s \in \sS 
\end{equation}

\begin{equation}\label[consset]{unmet_demand}
    \sum_{s \in \sS} x_{sr} + u_{r} = 1  \quad \forall r \in \sR
\end{equation}

\begin{equation}\label[consset]{flow_deviation3}
    g_{b} = \left | \frac{\sum_{(s,r)\in\sK_b} D_{r} x_{sr}}{\sum_{r\in\sR} D_r}  - Q_{b} \right |, \quad \forall b \in \sB
\end{equation}

\begin{equation}\label[consset]{non_negativity}
    x_{sr}, u_{r}, g_{b}  \in \mathbb{R}_{\geq 0}, x_{sr}, u_{r} \in [0,1]
\end{equation}

In objective function \labelcref{obj_fun_ss}, we minimize the weighted cost of i) unmet demand, ii) transportation cost for met demand, iii) supplier-receiver rating factor, and iv) percentage gap in demand tonnage in distance bins. Notice that, here, we condense variable $x_{sro}$ by removing $o$ index. Constraints \labelcref{Supplier}--\labelcref{flow_deviation3} function similar to the ones in \labelcref{FM_Supplier}--\labelcref{FM_flow_deviation} with the condensed variable $x_{sr}$. Constraints \labelcref{non_negativity} define variable domains.

\subsubsection{Commodity Assignment Model}

Now that we have solutions to $x_{sr}$ variables from the supplier selection model, we use them to assign commodity types to the pair of supplier and receiver $(s,r)$. The model is as follows.



\begin{equation}\label[objfunc]{obj_fun}
\begin{split}
    min_{y}  \quad w_5 \sum_{z \in \sZ, z^ \prime \in \sZ , o \in \sO} y_{z z^\prime o}
\end{split}
\end{equation}

\begin{equation}\label[consset]{unmet_commodity_demand}
    \sum_{o \in \sO} c_{sro} = 1,  \quad \forall s \in \sS, r \in \sR  
\end{equation}

\begin{equation}\label[consset]{FAF_deviation}
    y_{z z^\prime o} = \left | \sum_{\substack{s \in \sS|Z_s=z,\\ r \in \sR|Z_r=z^\prime}} W_{sr} c_{sro}  - P_{z z^\prime o}  \right |, \quad \forall z, z^ \prime \in \sZ, o \in \sO
\end{equation}

\begin{equation}\label[consset]{commodity_non_negativity}
    c_{sro}, y_{z z^\prime o}  \in \mathbb{R}_{\geq 0}, c_{sro} \in [0,1]
\end{equation}

Objective function \labelcref{obj_fun} minimizes the gap in the tonnage flow of commodities flowing between zones. Constraints \labelcref{unmet_commodity_demand} ensures integrality of assignments for commodity types to $(s,r)$ pairs. Constraints \labelcref{FAF_deviation} define $y_{zz^\prime o}$ in terms of $c_{sro}$ variables. Finally, \labelcref{commodity_non_negativity} define variable domains.

For the supplier rating, \citep{pourabdollahi_2017} proposed using the following proxies to substitute in the ordered logit model estimated based on the collected real-data: the unit value of the commodity from FAF data as a proxy for cost/price, production capacity as a proxy for capacity/reliability, and annual value of commodities as a proxy for credit/financial condition (which can be substituted with estimated revenues from the firm synthesizer). 
For the shipping cost, a travel distance and time matrix was estimated through POLARIS router module. In this paper, the cost was used as distance between shipper and receivers which accounts for network conditions, however, the cost can be also calculated based on routed travel time or a combination of both.

\begin{algorithm}
\caption{Supplier Selection and Commodity Assignment Problem Decomposition}\label{supplier_selection_algo}
\SetArgSty{textnormal}
\scriptsize
\SetKwInOut{Input}{Input}
\SetKwInOut{Output}{Output}
\Input{~ Establishment data, Make-Use table, FAF zonal flows, Shipping costs, Distance bin distribution.}
\Output{~ Supplier-receiver commodity flow assignments.}

\SetAlgoLined
\SetKwFunction{FFlowAssignment}{\textproc{FlowAssignment}}
\SetKwProg{Fn}{Function}{:}{}

\Fn{\FFlowAssignment{}}{
Initialize internal and external establishments: \\
    \quad $E_{\text{internal}} \gets$ internal establishments in study area  \\
    \quad $E_{\text{external}} \gets$ domestic external establishments \\
    \quad $R_{\text{internal}}, S_{\text{internal}} \gets$ internal receivers and suppliers  \\ 
    \quad $R_{\text{external}}, S_{\text{external}} \gets$ external receivers and suppliers 

\For{each receiver $r \in R_{\text{internal}} \cup R_{\text{external}}$}{
    $S_r \gets$ possible supplier sectors from Make-Use table  \tcp{$S_r$ is defined as sector of receiver.} 
    $Candidate\_Suppliers_r \gets \emptyset$ \\ 
    \For{each sector $\in Sector_r$}{
        Add suppliers in that sector and relevant FAF zones to $Candidate\_Suppliers_r$
        }
    }

\For{each pair $(r, s)$ where $s \in Candidate\_Suppliers_r$}{
    $C_{sr} \gets$ compute shipping cost between $s$ and $r$ \\
    $N_{sr} \gets$ estimate supplier rating \\
    $K_b \gets$ assign pair $(r, s)$ to a  distance bin $b$
}

Divide problem into subproblems: \\
    \quad $P_{II} \gets \{(s,r)\ |\ s \in S_{\text{internal}},\ r \in R_{\text{internal}}\}$ \\
    \quad $P_{EI} \gets \{(s,r)\ |\ s \in S_{\text{external}},\ r \in R_{\text{internal}}\}$ \\
    \quad $P_{IE} \gets \{(s,r)\ |\ s \in S_{\text{internal}},\ r \in R_{\text{external}}\}$

Solve Internal-Internal and External-Internal subproblems: \\
    \quad $Solution_{II} \gets$ solve supplier selection assignment over $P_{II}$  \\
    \quad $Solution_{EI} \gets$ solve supplier selection assignment over $P_{EI}$ \\

Update internal supplier production based on $Solution_{II}$ and $Solution_{EI}$

Solve Internal-External subproblem: \\
    \quad $Solution_{IE} \gets$ solve supplier selection assignment over $P_{IE}$ using updated supply

Combine all solutions: \\
    \quad $FinalSolution \gets Solution_{II} \cup Solution_{EI} \cup Solution_{IE}$

  \KwRet{$\mathit{flow\_assignment}$}
}
\end{algorithm}

\subsection{International Heuristic for Importer and Exporter Establishments Selection}

The major issue with supplier selection modeling in international trade lies in the lack of information on foreign establishments. Without their attributes, it is not feasible to solve the same supplier selection problem. Another complicating factor in international trade stems from the type of commodities, size, and business models of importers/exporters, which cannot be modeled easily given the data limitations. 

This probabilistic heuristic was developed specifically to address the unique challenge of linking aggregate port-level data with disaggregate firm-level agents in the absence of ground-truth micro-data. While it relies on established principles of flow disaggregation, the specific procedural logic is bespoke to the available data structure. A key advantage of this tailored approach is computational efficiency; as demonstrated in the results.

\cref{generate_import_export_algo} shows the heuristic used for selecting importer and exporter establishments for international shipments. The main inputs for this heuristic are (i) imported and exported tonnage by commodity type at each major port for each internal county, (ii) internal establishments for a given metropolitan area, (iii) production and consumption NAICS commodity matrices, (iv) size threshold or percentage of importer/exporter establishments, and (v) lower and upper bound of the trade volume for a given importer/exporter. The size threshold for establishments is based on the assumption in \cite{holguín-veras_2025} that mostly large sized establishments are involved in large volumes of international trade (neglecting international packages). Trade volume ranges are used as soft constraints on the quantity of goods imported/exported by a given establishment, to avoid over-allocation of goods to fewer establishments. The lower bound used in this study was a fully loaded truck once annually, while the upper bound was assumed to be four fully loaded trucks daily. These bounds are inputs to the algorithm, so they can be changed per the modeler discretion. These bounds act as soft constraints, since establishments who exceeds the range are not eliminated from the potential importer/exporter set to deal with instances of limited number of large sized businesses in a given county. Conversely, establishment trade volume by port can be below the range in case of low throughput ports. Assigning tonnage proportionally to the size of the importer/exporter is an alternative solution to this issue.

Once the importers/exporters set is defined, for both trade types (imports and exports) and for each commodity, a port is chosen out of the ports list. The NAICS of the importer/exporter is chosen based on the probability that this given NAICS produces/consumes this commodity. An establishment that belongs to the selected NAICS is chosen randomly and is assigned a flow within the specified range unless the port remaining tonnage is less than the lower bound. The shipment information is stored and the process is repeated till all the international trade volume has been assigned. 

\begin{algorithm}[!ht]
\caption{Import and Export Shipments Heuristic}\label[algo]{generate_import_export_algo}
\SetArgSty{textnormal}
\scriptsize
\SetKwInOut{Input}{Input}
\SetKwInOut{Output}{Output}
\Input{~ Ports, trade type, FAF import and export flows, commodity-industry sector mapping.}
\Output{~ Import and export commodity shipments between international ports and domestic establishments.}

\SetAlgoLined
\SetKwFunction{FGenerateShipments}{\textproc{GenerateImportExportShipments}}
\SetKwProg{Fn}{Function}{:}{}

\Fn{\FGenerateShipments{$\textit{trade\_type}$}}{

$importer\_exporters \gets$ filter establishments larger than a specified threshold; \\
$lb, ub \gets$ annual trade bounds;

\ForEach{$zone:  z \in \sZ$}{
 \ForEach{$commodity:  c \in \sC$}{
  $commodity\_production\_dict \gets$ industry sector production shares of $commodity\ c$ for $ zone \ z$\\
  $commodity\_consumption\_dict \gets$ industry sector consumption shares of $commodity\ c$ for $ zone \ z$ \\
  }
}
  
\ForEach{$trade\_type:  t \in \sT$}{
  \uIf{$trade\_type = \texttt{"export"}$}{
   $zone$ of establishment $\gets$ \texttt{origin}; \\
    $commodity\_dict \gets commodity\_production\_dict$; 
  }
  \uElseIf{$trade\_type = \texttt{"import"}$}{
    $zone$ of establishment $\gets$ \texttt{destination}; \\
    $commodity\_dict \gets commodity\_consumption\_dict$; 
  }

  \ForEach{$commodity:  c \in \sC$}{
    \texttt{ports\_flow$_{ct}$} $\gets$ flow of ports with $commodity = c$ and $trade\_type =  t$;

    \If{\texttt{ports\_flow$_{ct}$} $= \emptyset$}{
      \textbf{continue};
    }

    \ForEach{\texttt{port\_flow$_{ct}$} $\in$ \texttt{ports\_flow$_{ct}$}}{

      \texttt{industry\_sector\_shares} $\gets$ $commodity\_dict[c,$ \texttt{port\_flow$_{ct}[zone]$}];

      $importer\_exporters\_set \gets$ filter importers and exporter establishments where sector $\in$ sectors of \texttt{industry\_sector\_shares};

      \While{ \texttt{port\_flow$_{ct}[tons]$} $> 0$}{

        $selected\_sector \gets$ Industry sector randomly selected based on \texttt{industry\_sector\_shares};

        $sampled\_importer\_exporter \gets$ random sample from $importer\_exporters\_set$ matching $selected\_sector$ and $zone =  \texttt{port\_flow$_{ct}[zone]$}$;

        Sample $annual\_tons \sim \mathcal{U}(lb, ub)$;

        $annual\_demand \gets  \min($  \texttt{port\_flow$_{ct}[tons]$} $, annual\_tons)$ 
        
        \uIf{$trade\_type = \texttt{"export"}$}{
          $supplier \gets$ establishment ID from $sampled\_importer\_exporter$; \\
          $receiver \gets$ \texttt{port\_flow$_{ct}[international\_port]$}; 
        }
        \uElseIf{$trade\_type = \texttt{"import"}$}{
          $supplier \gets$ \texttt{port\_flow$_{ct}[international\_port]$}; \\
          $receiver \gets$ establishment ID from $sampled\_importer\_exporter$; 
        }

        Create $new\_shipment$ with $supplier$, $receiver$, $commodity$, and $annual\_demand$  \\

        Append $new\_shipment$ to $international\_shipments$;

        \texttt{port\_flow$_{ct}[tons]$} $\gets$ \texttt{port\_flow$_{ct}[tons]$} $- annual\_demand$;
      }
    }
  }
}
  \KwRet{$international\_shipments$}
}
\end{algorithm}

\section{Results and Discussion}\label{results_discuss}

\begin{figure*}[ht]
    \centering
    \includegraphics[width=\textwidth]{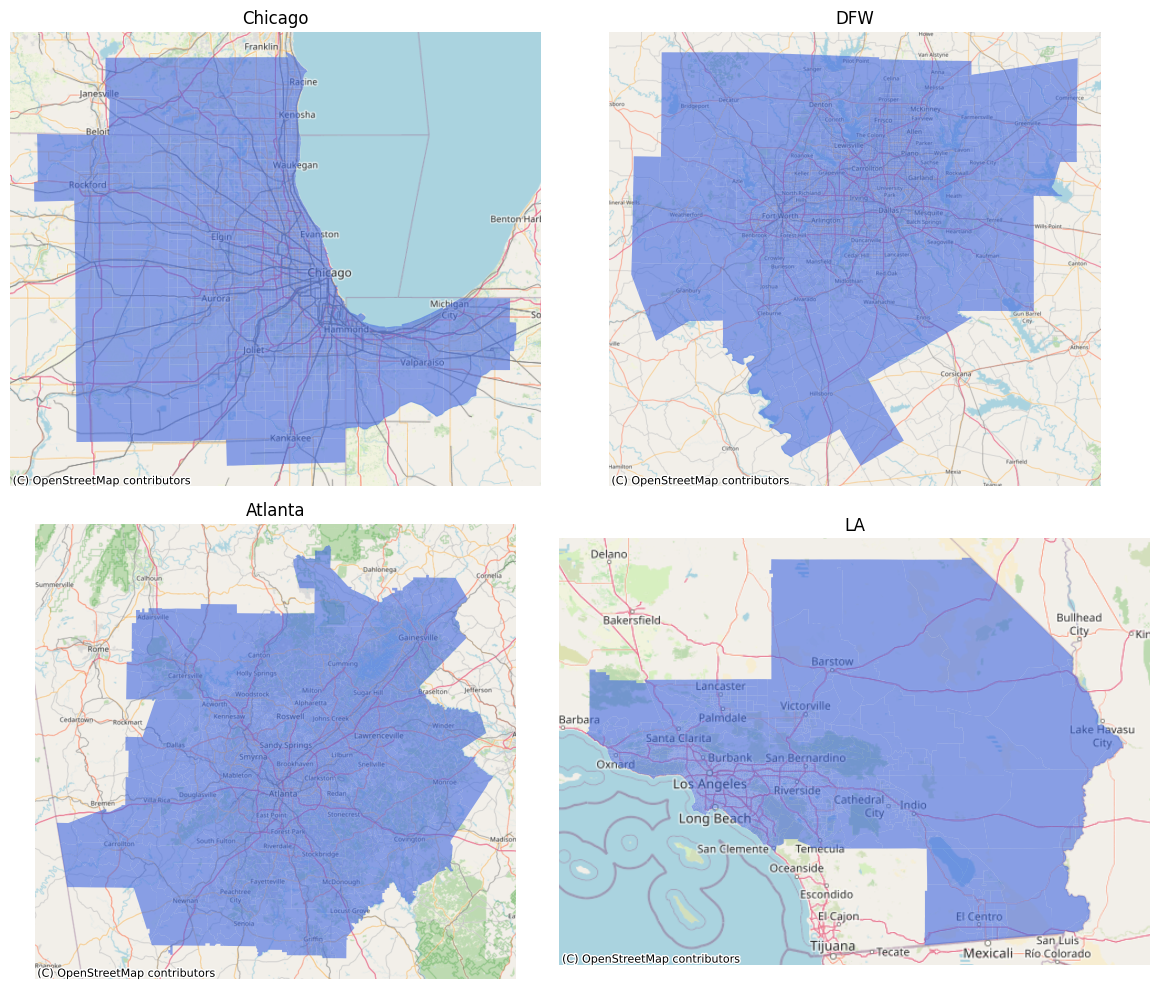}
    \caption{Modeled Metro Areas}
    \label[fig]{metro_area_maps}
\end{figure*}

The developed models were implemented in the metro areas of Atlanta, Chicago, DFW, and LA. \cref{metro_area_maps} shows the geographical regions of each study area, while \cref{runtime_summary} lists the run times for the used algorithms. It is important to provide a breakdown of the total computational effort. The times reported in \cref{runtime_summary} represent the core optimization and heuristic phases. However, the data preparation phase, which generates the set of potential supplier-receiver pairs and calculates their respective attribute matrices, remains computationally intensive. Even with the utilization of distributed multi-threading across multiple workstations, this pre-processing step accounted for approximately 60--70\% of the total end-to-end execution time. Consequently, a primary opportunity for future research lies in accelerating this data generation phase beyond standard CPU parallelization. Future efforts could leverage GPU-accelerated computing to handle these large-scale matrix operations or employ spatial indexing heuristics (e.g., KD-trees) to efficiently prune the candidate set prior to cost calculation, offering a greater marginal return on speed. All optimization-related computations were carried out on a single Intel{\textsuperscript \textregistered} Core\textsuperscript{TM} i9-14900K CPU @3.20 GHz workstation with 128 GB of RAM and 24 cores. Problem instances were solved by using
the Python 3.10.11 interface to the commercial solver Gurobi 11.0.3 \cite{gurobi}.

\cref{demand_summary} summarizes: (i) annual demand in million metric tons, (ii) international trade share, (iii) number of modeled internal and external establishments, (iv) final number of assignments based on 2023 FAF non-pipeline flows, (v) average number of domestic suppliers per receiver, (vi) average number of internal importers per port, and (vii) average number of internal exporters per port. \replace{}{To illustrate the computational burden, the total number of $x_{sr}$ variables evaluated during the optimization phase exceeded 393 million for Atlanta, 668 million for Chicago, 499 million for DFW, and 1 billion for LA.} It should be noted that since metro areas do not align perfectly with FAF zones, the FAF flows were disaggregated to reflect flow values of the study region counties. Atlanta has the lowest domestic and international demand, and lowest number of establishments and trade assignments. Chicago has the highest demand, however, LA has the most establishments and trade assignments. It is important to highlight that not all establishments are both goods producers and receivers, for example, hotels do not produce goods. Also, some of the external establishments send products to the study region, but not necessarily receive goods from the same region. The average number of domestic suppliers for a given receiver ranged between 1.17 to 1.32 for all metro areas. 

Atlanta and Chicago had the lowest averages for importers/exporters ratio, which is aligned with their lower share of international trade. Note that although Chicago generally has a very high import demand, this demand is mostly oil imported through pipelines which is not considered in this study due to the different nature of pipeline flows and their insignificance on the highway and railway networks. DFW had a high average number of exporters as it is a huge exporter metro area, while LA had a higher importer average for its high import volumes. The reported averages reflect the per-port values and are consequently affected by the presence of low-throughput ports. These results show that LA has a higher chance of getting impacted by import policies and disruptions, which would not only affect its position as a major gateway to the U.S., but also would impact the large number of businesses that strive on the import industry. Conversely, around 95\% of DFW's export tons are handled through points of entries within Texas itself, with almost 60\% of the tonnage and 80\% of exporters depending on Houston ports. These exporters are particularly vulnerable to disruptions during hurricane seasons, which have occasionally caused significant delays at the Port of Houston. This underlines the importance of modeling supplier selection problem to quantify impacts of possible disruptions on these freight flows and study mitigating measures. Moreover, modeling these trade partnerships allow for the assessment of potential policies and scenarios on the transportation network, e.g., the impact of increasing rail share for exports along the DFW-Houston corridor on the transportation network.  

\begin{table*}[ht]
\centering
\caption{Solver and Heuristic Runtimes by Metro Area}
\label[tab]{runtime_summary}
\begin{tabular}{l|c|c|c|c}
\hline
\multicolumn{1}{c|}{\multirow{2}{*}{\textbf{Metro Area}}} &
\multicolumn{2}{c|}{\textbf{Gurobi Solver Runtime (min.)}} &
\textbf{Imports/Exports} & \textbf{Total Runtime} \\
\cline{2-3} 
\multicolumn{1}{c|}{} & 
\textbf{Supplier Sel.} & 
\textbf{Comm. Assign.} & 
\textbf{Heuristic (min.)} & \textbf{(hrs)} \\
\hline
Atlanta  & 309.1 & 1.0 & 3.8 & 5.23 \\
Chicago  & 520.2 & 1.6 & 4.5 & 8.77 \\
DFW      & 406.8 & 1.1 & 4.9 & 6.88 \\
LA       & 807.7 & 1.4 & 2.5 & 13.53 \\
\hline
\end{tabular}
\end{table*}

\begin{table*}[ht]
\centering
\caption{Metro Area Key Statistics}
\label[tab]{demand_summary}
\begin{tabular}{l|c|c|c|c|c|c}
\hline
\multicolumn{1}{c|}{\multirow{2}{*}{\textbf{Metro Area}}} & 
\multicolumn{1}{c|}{\textbf{Annual $10^6$ tonnage}} & 
\multicolumn{1}{c|}{\multirow{2}{*}{\textbf{Establishments}}} & 
\multicolumn{1}{c|}{\textbf{Trade}} &
\multicolumn{3}{c}{\textbf{Average}} \\
\cline{5-7} 
\multicolumn{1}{c|}{} & 
\textbf{(International \%)} & 
\multicolumn{1}{c|}{} & 
\textbf{Assignments} &
\textbf{Sup.} & 
\textbf{Imp.} & 
\textbf{Exp.} \\
\hline
Atlanta  & 231.8 (10.3\%)  & 144{,}705 & 267{,}010 & 1.2 & 39 & 48\\
Chicago  & 539.3 (8.3\%) & 229{,}136 & 424{,}334 & 1.3 & 40 & 52 \\
DFW      & 466.8 (18.4\%) & 169{,}409 & 342{,}488 & 1.3 & 124 & 492 \\
LA       & 423.3 (16.8\%) & 326{,}587 & 618{,}515 & 1.2 & 417 & 38 \\
\hline
\end{tabular}
\end{table*}

\begin{figure*}[!ht]
    \centering
    \includegraphics[width=\linewidth]{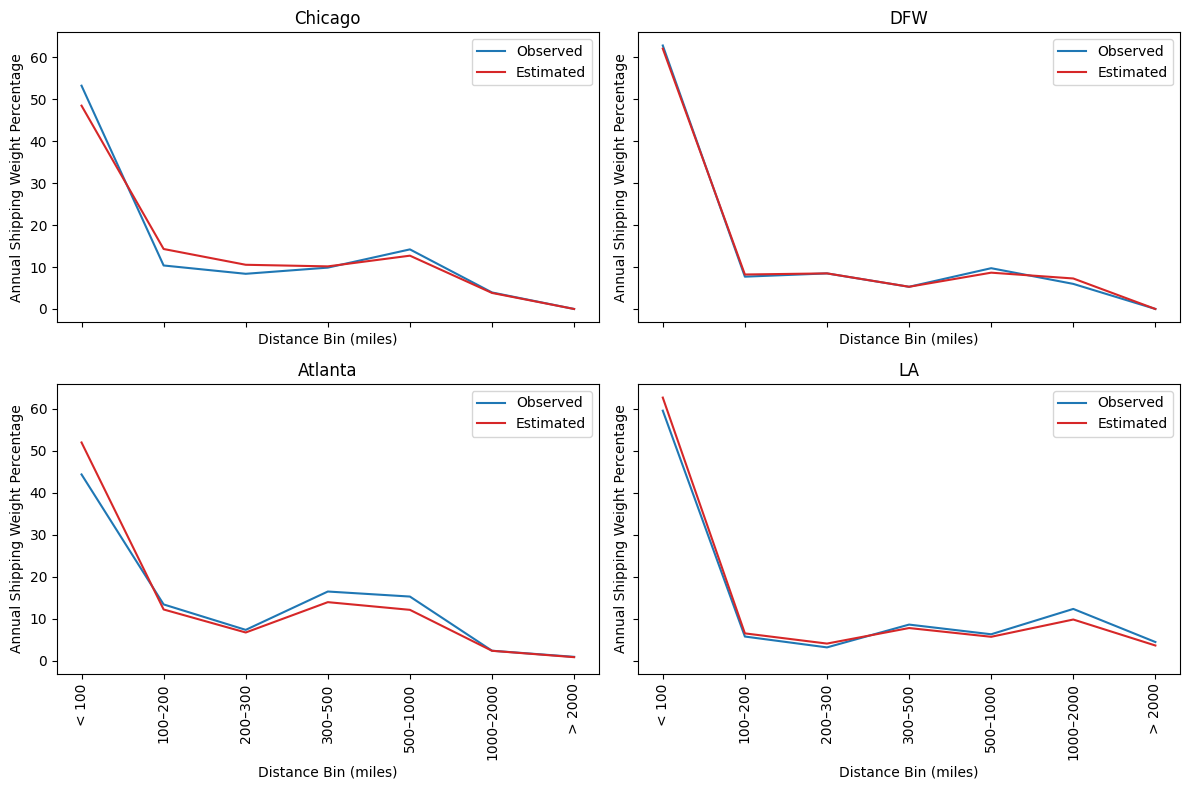}
    \caption{Shipping Distance Distribution}   
    \label[fig]{ship_dist}
\end{figure*}

\cref{ship_dist} illustrates the regional domestic shipping distance distribution for both observed and optimized flows. The observed and estimated weight percentages are recorded in \cref{distance_bin_comparison}, where $\Delta$ refers to the (observed - estimated) percentages. \replace{}{In this implementation, the same aggregate distance bins and target distributions are applied universally across all industry types. While industry-specific distance distributions could theoretically be utilized to capture sectoral heterogeneity, the sample size of regional CFS data is often sparse to reliably estimate such distributions at the disaggregated NAICS level per each region.} The highest $\Delta$ percentage happened in Atlanta where the model overestimated by 7.6\% in the first distance bin, the rest of the bins in all the cities have $\Delta$ less than 5\%. \replace{}{To quantify the aggregate goodness-of-fit, the Mean Absolute Error (MAE) relative to the percentage shares across the seven bins was calculated, yielding 2.2\% for Atlanta, 1.8\% for Chicago, 0.5\% for DFW, and 1.4\% for LA. For instance, the MAE for Chicago is obtained by averaging the absolute percentage errors ($\Delta$) from its seven distance bins: $(|4.8| + |-3.9| + |-2.2| + |-0.3| + |1.5| + |0.1| + |0.0|) / 7 \approx 1.8\%$. It is important to note that this calibration is achieved at the aggregate regional level; specific commodity groups may exhibit spatial heterogeneity, meaning the distance distributions for specialized goods (e.g., petroleum) might diverge from this aggregate trend.} Overall, 40–60\% of the annual shipments occur within distances of less than 100 miles, aligning with the national-level freight patterns. This result is also consistent with FAF statistics, where intra-zonal flows less than 100 miles account for approximately 60\% of the total national freight tonnage in U.S. metropolitan areas. This demonstrates that freight shipments in major U.S. metropolitan areas are predominantly short-distance, with a sharp decline as distance increases. This trend reinforces the need for policies targeting efficient local and regional freight movement. 

It should be noted that using smaller bins, slightly increases the gaps, however, the major increase occurs in bins below 50 miles. This is due to using establishments' county centroids to calculate distance in this model and then comparing the results to actual observed distances from the CFS survey, which although obscures the actual location of the supplier and receiver, reports the observed distances. Such differences might not impact the results of longer shipments but significantly affects shorter shipments especially inter-county shipments. A possible solution to this issue is to use POLARIS disaggregated transportation analysis zones or locations, which can help reduce bias arising from large inter-county distances. Yet, improving the resolution of the model brings in further computational challenges contradicting with the simplification approach aimed to be followed in this study. It is worth mentioning that the gaps cannot be entirely eliminated since the model has multi-objectives and the priority was set to ensuring receiver demands are met. Also, the external establishments sampled in the FG module in POLARIS does not consider the shipping distances, introducing an initial bias to the model.

\begin{table*}[ht]
\centering
\caption{Observed vs. Estimated Distance Bin Percentage Shares}
\label[tab]{distance_bin_comparison}
\begin{tabular}{l|ccc|ccc|ccc|ccc}
\hline
\multirow{2}{*}{\begin{tabular}[c]{@{}c@{}}\textbf{Distance}\\ \textbf{bin (mi)} \end{tabular}} &
\multicolumn{3}{c|}{\textbf{Atlanta}} &
\multicolumn{3}{c|}{\textbf{Chicago}} &
\multicolumn{3}{c|}{\textbf{DFW}} &
\multicolumn{3}{c}{\textbf{LA}} \\
\cline{2-13}
 & \textbf{Obs.} & \textbf{Est.} & \textbf{$\Delta$} 
 & \textbf{Obs.} & \textbf{Est.} & \textbf{$\Delta$}
 & \textbf{Obs.} & \textbf{Est.} & \textbf{$\Delta$}
 & \textbf{Obs.} & \textbf{Est.} & \textbf{$\Delta$} \\
\hline
$< 100$      & 44.4 & 52.0 & -7.6  & 53.2 & 48.5 & 4.8   & 62.8 & 62.1 & 0.7  & 59.6 & 62.7 & -3.1 \\
100--200     & 13.4 & 12.2 & 1.2   & 10.4 & 14.3 & -3.9  & 7.7  & 8.2  & -0.5   & 5.7  & 6.5  & -0.8 \\
200--300     & 7.3  & 6.7  & 0.6   & 8.4  & 10.5 & -2.2  & 8.5  & 8.5  & -0.01  & 3.1  & 4.0  & -0.9 \\
300--500     & 16.5 & 13.9 & 2.5   & 9.9  & 10.2 & -0.3  & 5.3  & 5.3  & 0.00  & 8.6  & 7.8  & 0.8 \\
500--1000    & 15.3 & 12.1 & 3.2   & 14.2 & 12.7 & 1.5   & 9.7  & 8.7  & 1.1  & 6.3  & 5.7  & 0.6 \\
1000--2000   & 2.3  & 2.3  & 0.01  & 3.9  & 3.8  & 0.1   & 6.0  & 7.3  & -1.3   & 12.3 & 9.8  & 2.5 \\
$> 2000$     & 0.88 & 0.81 & 0.08  & 0.002 & 0.002 & 0.0 & 0.01 & 0.00 & 0.01 & 4.42 & 3.59 & 0.83 \\
\hline
\end{tabular}
\end{table*}

\begin{figure*}[!ht]
    \centering
    \includegraphics[width=\linewidth]{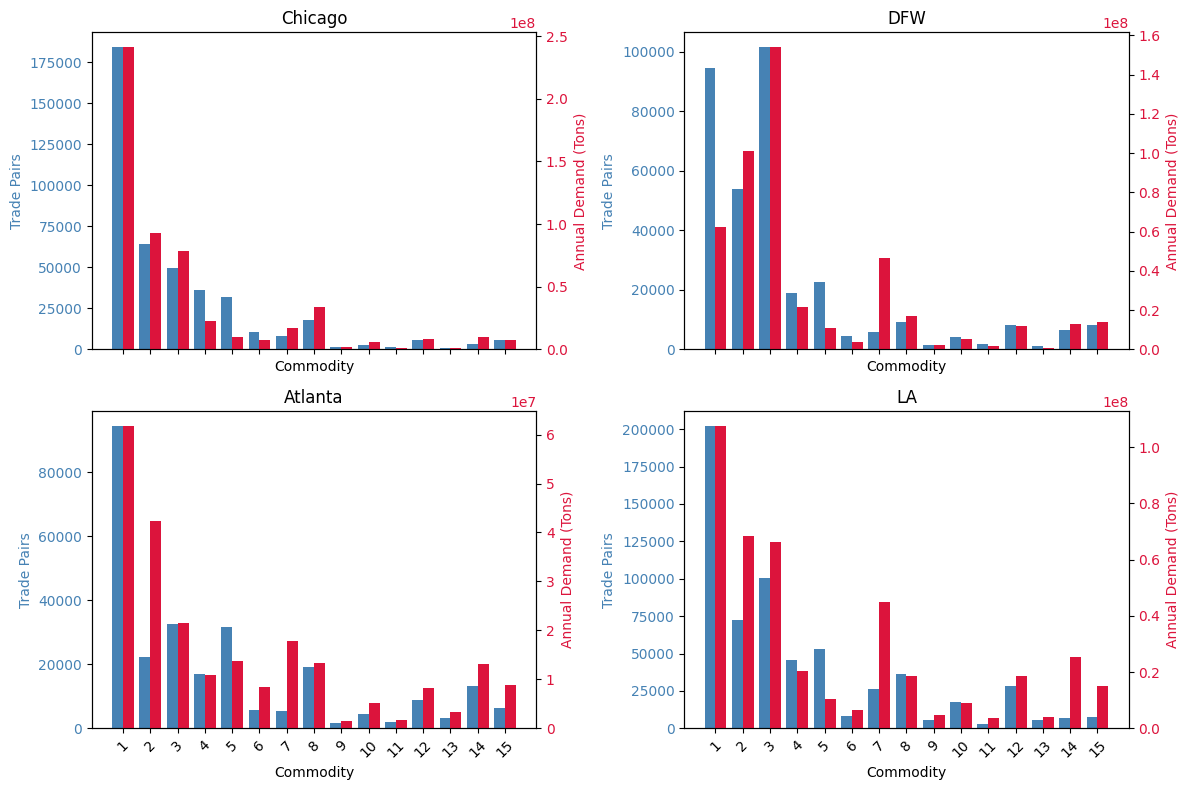}
    \caption{Trade Assignments and Demand by Commodity Type}   
    \label[fig]{commodities}
\end{figure*}

\cref{commodities} depicts the analysis of number of trade assignments and annual demand by commodity grouping (shown in \cref{commodity_groups}), revealing significant spatial sector variations among major U.S. metropolitan regions. In Atlanta, Chicago, and LA, food and agricultural commodities—including processed foods—consistently represent both the largest share of trade pairs and the highest total demand. This dominance underscores the centrality of these commodities to urban consumption patterns, particularly given their primary end-users in retail, restaurants, and food service sectors. DFW diverges from this trend, with petroleum products surpassing food commodities in both trade volume and aggregate demand. This reflects the DFW’s industrial structure and the major role of the petroleum sector in shaping freight flows.

Population size is the key driver of total food commodity demand. Chicago and LA, the most populous metro areas in this study, exhibit the highest aggregate demand for food products. The number of businesses of NAICS 722-Food Services and Drinking Places-can be used as a proxy to analyze the dominance of the food commodity group. For instance, Chicago, with the highest demand tonnage and approximately 19,000 food and drinking places, averages 1,310 annual tons per trade pair. In contrast, LA, despite a high total demand, has nearly 34,000 such establishments, resulting in a lower average of 532 tons per pair. Atlanta and DFW, with 11,000 and 14,000 establishments respectively, display intermediate averages 654 and 662 tons per pair, reflecting both their relatively smaller populations and lower total demand. These findings suggest that while demand for food commodities scales with metropolitan population, the average demand per trade pair is inversely related to the number of establishments due to higher market competition.

Finally, commodities such as leather, textiles, electronics, and office furniture consistently shows the lowest trade volumes across all regions, indicating their relatively minor contribution in terms of volume, and consequently in terms of freight flows on the transportation network. These results show that a detailed understanding of commodity-specific trade dynamics provides a better modeling tool suited to study the appropriate targeting policies for a given metropolitan area.

\section{Conclusion}\label{conclusion}

The freight transportation field often suffers from data limitations, which constrain modeling capabilities, making researchers rely on publicly available aggregated data to infer the actual behaviors and operations of freight stakeholders. One such challenge arises in the supplier selection problem, where data on actual bidding processes and firm-level decision-making are rarely accessible, making it difficult to accurately represent market conditions. This research addresses the supplier selection process by modeling receiver behavior to maximize perceived supplier ratings, but more importantly by decreasing the discrepancy between modeled and observed commodity flows. The overarching objective is to produce a more behaviorally realistic and transport-sensitive representation of supplier selection outcomes, particularly in terms of their impact on transportation network flows. To achieve this goal, this paper has proposed supplier selection and commodity assignment model that seek to match shipping distance distribution while ensuring a match in the inter-zonal commodity flows. In addition to an international shipments heuristic that matches commodity flows individual ports with establishments. The developed models were implemented on a large scale on four metropolitan areas in the U.S.: Atlanta, Chicago, DFW, and LA. The model results showed a close match between the estimated and observed shipping distance distributions for all study regions.  

The implications of the proposed models are far-reaching, providing a critical tool for studying targeted policies and understanding complex freight dynamics. By capturing the nuanced decisions of individual freight agents, the model enables a deeper exploration of how changes in sourcing patterns, influenced by factors like trade policies, infrastructure investments, or even disruptions, translate into real-world transportation network impacts. For instance, the model can be utilized to assess the network-level consequences of cost fluctuations on specific imported goods by simulating shifts in supplier selection, evaluating changes in VMT, and quantifying their effects on congestion and network reliability. Conversely, it can inform strategic infrastructure planning by predicting how new or improved corridors might alter logistics choices, attract new suppliers, or facilitate more efficient movement of goods. Furthermore, the ability to simulate supplier-receiver relationships at a micro-level within POLARIS ABM makes this model invaluable for analyzing supply chain resilience strategies. It can help identify vulnerabilities within existing supply chains, simulate the cascading effects of disruptions (e.g., port closures, infrastructure failures), and evaluate the effectiveness of mitigation measures such as diversifying supply sources or rerouting commodity flows. This granular insight into trade partnerships offers policymakers a powerful analytical framework to proactively address emerging challenges and optimize the performance of urban and national freight systems.

A critical distinction must be made regarding the scalability of the joint versus decomposed formulations. In the joint formulation, the number of decision variables ($x_{sro}$) scales linearly with the number of commodity types ($|\mathcal{O}|$). For large metropolitan datasets containing millions of potential supplier-receiver pairs, this linear increase results in a problem matrix size that exceeds the memory capacity of standard high-performance workstations, rendering the joint model computationally intractable for detailed commodity disaggregation. 
Conversely, the proposed decomposed formulation creates a scalable pathway. By aggregating flows into $x_{sr}$ during the computationally intensive supplier selection phase, the bottleneck step becomes independent of the commodity set size. This ensures that the model's run-time is driven primarily by the density of the trade network (the number of feasible supplier-receiver pairs) rather than the granularity of the commodity classification.

However, for a unified national-scale implementation where the network size increases dramatically, the exact decomposition approach may still face tractability limits. Addressing this national scalability will likely require future research into heuristic solution methods or hierarchical regional decomposition strategies to trade off a marginal degree of optimality for necessary computational speed.

Despite its significant contributions, this research acknowledges several limitations that offer clear avenues for future work. A key consideration is that while the commodity assignment problem is made feasible through decomposition, the current approach does not explicitly guarantee an optimal commodity assignment from a global perspective, given its sequential nature after supplier selection. Future efforts could explore iterative or more integrated solution methodologies to enhance optimality. Another limitation stems from the spatial aggregation of shipping costs; currently, distances between establishments are based on county centroids, which, especially in large or geographically diverse counties, may not accurately reflect actual travel distances. A refined approach would involve leveraging the more precise POLARIS locations and zones, or even exact establishment coordinates, to calculate transportation costs, thereby increasing the model's accuracy. Additionally, the initial set of external establishments generated by the FG module in POLARIS introduces a potential bias, as these are sampled without explicit consideration of their optimal shipping distances. Future enhancements to the FG module could integrate spatial optimization to yield a more representative initial population of external trade partners. The present shipping cost calculations do not explicitly account for network congestion. However, a major future direction involves integrating the model with the POLARIS Freight multimodal router, which is currently under development. This integration would allow for dynamic feedback, where simulated congestion levels on the network would influence shipping costs, leading to more realistic and adaptive supplier selection decisions in subsequent simulation iterations. \replace{Finally, using the heuristic for imports and exports, although warranted due to the impracticality of getting data on foreign suppliers and receivers, does add a limitation to the work. As external assumptions must be made to model scenarios replicating unseen changes in foreign networks. However, the heuristic can still capture the domestic impacts of port disruptions.}{Finally, using the heuristic for imports and exports, although warranted due to the impracticality of getting data on foreign suppliers and receivers, introduces a limitation regarding parameter sensitivity. The resulting allocation of international flows is inherently sensitive to the assumed lower and upper trade bounds; tightening these bounds forces the flow across a larger pool of establishments, while relaxing them concentrates trade among a few dominant actors. Despite this sensitivity in baseline allocation, the heuristic remains a robust tool for capturing the relative domestic impacts of supply shocks, such as port disruptions.} Addressing these limitations will further enhance the model's fidelity, computational efficiency, and practical utility for advanced freight planning and policy analysis.

\section*{Acknowledgments}
This material is based on work supported by the U.S. Department of Energy, Office of Science, under contract number DE-AC02-06CH11357. This report and the work described were sponsored by the U.S. Department of Energy (DOE) Transportation Technologies Office (TTO) under the Integrated Transportation and Energy Cross-Sectoral System of Systems at Scale (ITES4), an initiative of the Energy Efficient Mobility Systems (EEMS) Program. Melissa Rossi, a DOE Office of Energy Critical Minerals and Energy Innovation (CMEI) manager, played an important role in establishing the project concept, advancing implementation, and providing guidance.

\section*{Conflict of Interest}
The authors declared no potential conflicts of interest with respect to the research, authorship, and/or publication of this article.

\section{Author Contributions}
The authors confirm their dual responsibility for the following: study conception and design, data collection, analysis and interpretation of results, and manuscript preparation.

\section*{ORCID}

\def\orcid#1{\kern .08em\href{https://orcid.org/#1}{\includegraphics[keepaspectratio,width=0.7em]{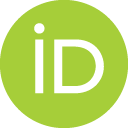}}}
Abdelrahman Ismael \orcid{0000-0003-0371-2110}

\noindent Taner Cokyasar \orcid{0000-0001-9687-6725}

\newpage

\clearpage
\bibliographystyle{TRR}
\bibliography{trb_template}

\vfill
\framebox{\parbox{.90\linewidth}{\scriptsize The submitted manuscript has been created by
        UChicago Argonne, LLC, Operator of Argonne National Laboratory (``Argonne'').
        Argonne, a U.S.\ Department of Energy Office of Science laboratory, is operated
        under Contract No.\ DE-AC02-06CH11357.  The U.S.\ Government retains for itself,
        and others acting on its behalf, a paid-up nonexclusive, irrevocable worldwide
        license in said article to reproduce, prepare derivative works, distribute
        copies to the public, and perform publicly and display publicly, by or on
        behalf of the Government.  The Department of Energy will provide public access
        to these results of federally sponsored research in accordance with the DOE
        Public Access Plan \url{http://energy.gov/downloads/doe-public-access-plan}.}}
\end{document}